\documentclass[12pt,thmsa]{article}
\usepackage{amsmath, latexsym, amsfonts, amssymb, amsthm, amscd}
\usepackage{graphics}

\newtheorem{theorem}{Theorem}
\newtheorem{corollary}{Corollary}

\newtheorem{lemma}{Lemma}

\newtheorem{rem}{Remark}

\newcommand{\p}{\Bbb{P}}

\newcommand{\e}{\Bbb{E}}

\newcommand{\ind}{\mbox{\rm 1\hspace{-0.04in}I}}

\newcommand{\ed}{\stackrel{(d)}{=}}

\begin{document}

\title{Invariance principles for local times at the supremum of random walks
and L\'{e}vy processes.}
\author{}
\maketitle

\begin{center}
{\large L. Chaumont}\footnote{
LAREMA, D\'{e}partement de Math\'{e}matiques, Universit\'{e} d'Angers, 2, Bd
Lavoisier - 49045, \newline
\hspace*{0.25in}\textsc{Angers Cedex 01.}$\;\;\;$E-mail:
loic.chaumont@univ-angers.fr$\;\;\;$} {\large and R.A. Doney}\footnote{
Department of Mathematics, University of Manchester, \textsc{Manchester,
\textrm{M 13 9 PL.}}\newline
\hspace*{.21in} E-mail: rad@ma.man.ac.uk}
\end{center}

\vspace{0.03in}

\begin{abstract}
We prove that when a sequence of L\'{e}vy processes $X^{(n)}$ or a normed
sequence of random walks $S^{(n)}$ converges a.s. on the Skorokhod space
toward a L\'{e}vy process $X$, the sequence $L^{(n)}$ of local times at the
supremum of $X^{(n)}$ converges uniformly on compact sets in probability
toward the local time at the supremum of $X$. A consequence of this result
is that the sequence of (quadrivariate) ladder processes (both ascending and
descending) converges jointly in law towards the ladder processes of $X$. As
an application, we show that in general, the sequence $S^{(n)}$ conditioned
to stay positive converges weakly, jointly with its local time at the future minimum,
towards the corresponding functional for the limiting process $X$. From this 
we deduce an invariance principle for the meander which extends known results for 
the case of attraction to a stable law.\newline

\newline

\end{abstract}

\vspace{0.1in}

\section{Introduction}

It is well-known that if a sequence of L\'{e}vy processes $X^{(n)}$
converges a.s. on the Skorokhod space to a limiting L\'{e}vy process $X$,
then the corresponding sequence of local times at a fixed level of $X^{(n)}$
do not necessarily converge to the local time of $X$, whatever the
definition the local times of $X^{(n)}$ is: occupation time, crossing
times,... However, in fluctuation theory of L\'{e}vy processes, it is the
local times at extrema which play a major r\^{o}le, not the local times at
fixed levels, so a natural and important question is whether these local
times converge. A similar question can be posed about the local times at
extrema of a sequence of normed random walks which converge to \ a L\'{e}vy
process.

To our knowledge, the only results known in this vein can be found in
Greenwood, Omey and Teugels\cite{got} and in Duquesne and Le Gall \cite{dlg}.
The first paper deals with the "classical" case where $S^{(n)}$ is got by
norming a fixed random walk $S,$ the assumption being that for some norming
sequence $c_{n,}$ $(S_{[nt]}/c_{n},t\geq 0)$ converges in law to $X,$
necessarily stable, and the conclusion is that a normed version of the
bivariate ladder process of $S$ converges in law to the bivariate ladder
process of $X$. One easily derives that a normed version of the local time
at the maximum of $S$ converges in law to the local time at the supremum of $
X.$ (A different proof of this result and a converse result can be found in
Doney and Greenwood \cite{dg}.)

The second paper considers a more general scenario where each $S^{(n)}$ is
got by norming a different random walk, but restricts itself to the case
that each random walk is \ downwards skip-free, so that the limiting L\'{e}vy
process is automatically spectrally positive. (This is because the
result, Theorem 2.21 of \cite{dlg}, is a tool for the study of the height
process of the sequence of Galton-Watson processes related to $S^{(n)}$.)
Again convergence in law is assumed, and the conclusion is again convergence
in law of a normed version of the local time.

In this paper, we give three major extensions of these results. In Theorem 1
we show that whenever $S^{(n)}$ converges in law to $X,$ then a normed
version of the bivariate ladder process of $S$ converges in law to the
bivariate ladder process of $X$. Again, we can deduce that a normed version
of the local time at the maximum of $S$ converges in law to the local time
at the supremum of $X.$ (Our only assumption on $X$ is that it has a
continuous local time $L$ at the supremum, but if this were to fail a
similar result could be formulated.) Next, in Theorem 2, we show that if the
assumption is strengthened to a.s. convergence, then the normed sequence of
local times converge to $L$ in probability uniformly on compacts. This
result allows us to deduce, in Theorem 3, an analogous result when a
sequence $X^{(n)}$ of L\'{e}vy processes converges a.s. to $X.$ (We stress
that for such a result to hold, we have to remove the ambiguity inherent in
the definition of local times for L\'{e}vy processes by insisting on a
standard normalisation for the local times of $X^{(n)}$ and $X$: see (\ref
{norm}).) An important corollary of these results is the convergence in law
of the quadrivariate process of upgoing and downgoing ladder processes: see
Corollary \ref{Coro4}.

In the last section, we show that if a sequence $(S^{(n)}_{[nt]},t\ge0)$ of
continuous time random walks converges in law toward a L\'evy process $X$,
then the sequence of these processes  conditioned to stay positive on the
whole time interval $[0,\infty)$ converges in law toward $X$ conditioned to
stay positive. We illustrate the usefulness of the results of Section \ref
{main} by showing, that this convergence also holds jointly with the local
time at the future minimum. Finally, we obtain an invariance principle for
the meander, i.e. we show that the sequence $(S^{(n)}_{[nt]},t\ge0)$
conditioned to stay positive over $[0,1]$ converges in law towards $X$
conditioned to stay positive over $[0,1]$. These results extend the
"classical" case studied by Bolthausen \cite{bo}, Doney \cite{do0} and
Chaumont and Caravenna~\cite{cc}.

\section{Preliminaries}

\label{prelim}

Let $X$ be any L\'{e}vy process for which $0$ is regular for the open
half-line $(0,\infty ).$ Then $0$ is also regular for itself for the
reflected process $R:=M-X$, where $M_{t}=\sup_{0\leq s\leq t}X_{s},$ and so
there exists a continuous local time for $R$ at $0.$ This local time $L$ is
only specified up to multiplication by a constant, but we will assume
throughout that its normalization is fixed by the requirement that
\begin{equation}
\mathbb{E}\left( \int_{0}^{\infty }e^{-t}\,dL_{t}\right) =1.  \label{norm}
\end{equation}
The process $L$ will be called the \textit{local time of $X$ at its maximum}.
It satisfies $L_\infty<\infty$, a.s. if and only if $X$ drifts to $-\infty$.
\newline

Let us introduce the ascending bivariate ladder process $(\tau ,H)$: the
ladder time process is $\tau _{t}=\inf\{s:L_s>t\}$, with the convention that
$\inf\emptyset=+\infty$ and the ladder height process is $H_{t}=X(\tau_{t})$,
if $\tau_t<\infty$ and $H_t=\infty$, if $\tau_t=\infty$. The process $
(\tau ,H)$ is a (possibly killed) bivariate subordinator whose Laplace
exponent is given by Fridstedt's formula:
\begin{equation*}
\kappa (\alpha ,\beta )=-\log \{e^{-(\alpha \tau _{1}+\beta H_{1})}\}=\exp
\left( \int_{0}^{\infty }\int_{0}^{\infty }(e^{-t}-e^{-\alpha t-\beta
x})t^{-1}\mathbb{P}\{X_{t}\in dx\}dt\right),
\end{equation*}
for $\alpha,\beta\ge0$, with the convention that $e^{-\infty}=0$. See Chapter VI of \cite{be} or
Chapter 4 of \cite{do}: note that (\ref{norm}) squares with $\kappa (1,0)=1.$
We write $q_H$, $\delta _{H}$ and $\pi^{H}$ respectively for the killing
rate, the drift coefficient and the L\'{e}vy measure of $H$. In particular,
the Laplace exponent of $H$ is given by
\begin{equation*}
\kappa(0,\beta)=q_H+\delta_H\beta+\int_0^\infty(1-e^{-\beta x})\,\pi^H(dx)\,.
\end{equation*}
Note that our assumptions imply that if $\delta_{H}=0$, then $
\pi^{H}(0,\infty)=\infty$, since otherwise $0$ would be irregular for the
open half-line $(0,\infty )$ for $X$.\newline

A random walk is a discrete time process $S=(S_{k},k=0,1,\cdots )$ such that
$S_{0}=0$ and, for $k\geq 1,$ $S_{k}=\sum_{1}^{k}Y_{r},$ where $
Y_{1},Y_{2},\cdots $ are independent and identically distributed. We define
the local time at its maximum of any random walk $S$ by $\Lambda_0=0$ and
for all $k\ge1$,
\begin{equation}
\Lambda _{k}=\#\{j\in \{1,\dots ,k\}:S_{j-1}<S_j,\,S_{j}=\max_{i\leq
j}S_{i}\}\,.  \label{loctimerw}
\end{equation}
As in continuous time, $\Lambda_\infty<\infty$, a.s. if and only if $S$
drifts to $-\infty$. We also introduce the strict ascending ladder processes
for $S$. The strict ascending ladder time process $T$ of $S$ is defined by $
T_{0}=0$ and for all $k\ge0$,
\begin{equation*}
T_{k+1}=\min \{j>T_{k}:S_{j}>S_{T_{k}}\}\,,
\end{equation*}
with $\min\emptyset=\infty$. The strict ascending ladder height process is
given by
\begin{equation*}
H_{k}=S(T_{k})\,,\;\;\mbox{if}\;\; T_k<\infty\;\;\mbox{and}
\;\;H_k=\infty\,,\;\;\mbox{if}\;\;T_k=\infty\,.
\end{equation*}
Note that $T$ is the inverse of $\Lambda$, i.e. $\Lambda_{T_k}=k$, for all $
k\le\Lambda_\infty$. Let us mention that all the results of this paper are
still valid if in the statements one replaces the \textit{strict} ladder
process and the \textit{strict} local time by the \textit{weak} ladder
process and the \textit{weak} local time.\newline

In the next sections, $S^{(n)}$ will denote a random walk whose distribution
can depend on $n$ and $\Lambda ^{(n)},T^{(n)}$ and $H^{(n)}$ will denote the
corresponding local time, ladder time, and ladder height process. We will
say that the sequence of random walks $S^{(n)}$ converges weakly (resp.
almost surely) toward the L\'{e}vy process $X$ if the sequence of continuous
time processes $(S_{[nt]}^{(n)},\,t\geq 0)$ converges weakly (resp. almost
surely) toward $X$ on the Sorokhod space $\mathcal{D}([0,\infty))$ of
c\`{a}dl\`{a}g paths. Note that according to Theorem 2.7 of Skorokhod \cite{sk},
if the process $(S_{[nt]}^{(n)},\,t\geq 0)$ converges in the sense of finite
dimensional distributions, then it converges weakly.
If a stochastic process $Y$ has lifetime $\zeta$ and if
the $Y^{(n)}$'s have lifetimes $\zeta^{(n)}$, then we say that the sequence $
Y^{(n)}$ converges toward $Y$ in some sense if the sequence of processes $
(Y_t^{(n)}\mbox{\rm 1\hspace{-0.04in}I}_{\{t<\zeta^{(n)}\}}+Y_{
\zeta^{(n)}-}^{(n)}\mbox{\rm 1\hspace{-0.04in}I}_{\{t\ge
\zeta^{(n)}\}},\,t\ge0)$ converges toward the process $(Y_t\mbox{\rm
1\hspace{-0.04in}I}_{\{t<\zeta\}}+Y_{\zeta-}\mbox{\rm 1\hspace{-0.04in}I}
_{\{t\ge \zeta\}},\,t\ge0)$ in this sense on the space $\mathcal{D}
([0,\infty))$. Note also that weak (resp. almost sure) convergence of stochastic
processes on the
space $\mathcal{D}([0,\infty))$ is equivalent to weak (resp. almost sure)
convergence on the space $\mathcal{D}([0,t])$ for all $t>0$, see Theorem
16.7 in \cite{bi}. Weak or almost sure convergence of a sequence of stochastic
processes $Y^{(n)}$ toward $Y$ will be denoted respectively by $Y^{(n)}
\overset{(\mbox{\rm\tiny law})}{\longrightarrow }Y$ and $Y^{(n)}\overset{(
\mbox{\rm\tiny a.s.})}{\longrightarrow }Y$.

\section{Main results}

\label{main}

The following result extends Theorem 3.2 in \cite{got} and lemme 3.4.2, p.54
of \cite{vi}.

\begin{theorem}
\label{Th1} Let $X$ be any L\'{e}vy process such that $0$ is regular for the
open half line $(0,\infty )$ and assume that some sequence of random walks $
S^{(n)}$ converges in law toward $X$. Then we have the following convergence
in law:
\begin{equation*}
\left[ \left( n^{-1}T_{[a_{n}t]}^{(n)},H_{[a_{n}t]}^{(n)}\right) ,\,t\geq 0
\right] \overset{(\mbox{\rm\tiny law})}{\longrightarrow }(\tau ,H)\,,
\end{equation*}
as $n\rightarrow\infty$, where
\begin{equation}
a_{n}=\exp \left( \sum_{k=1}^{\infty }\frac{1}{k}e^{-k/n}\mathbb{P}
(S_{k}^{(n)}>0)\right).  \label{a}
\end{equation}
\end{theorem}

\begin{rem}
\label{rem1} Under the hypothesis of Theorem $\ref{Th1}$, i.e. when $0$ is
regular for $(0,\infty)$, Rogozin's criterion asserts that $\int_0^1t^{-1}
\mathbb{P}(X_t>0)\,dt=\infty$, see $\cite{be}$, Proposition VI.$3.11$. It
follows from this result and weak convergence of $S^{(n)}$ toward $X$ that
in Theorem $\ref{Th1}$, we necessarily have $\lim_{n\rightarrow\infty}a_n=
\infty$.
\end{rem}

\begin{rem}
The sequence $S^{(n)}$ could also be written in the form
\begin{equation*}
S^{(n)}=\frac{1}{c_{n}}\tilde{S}^{(n)}
\end{equation*}
and then we would recover the standard formulation for triangular arrays.
But in this case, using obvious notations, the result of Theorem $\ref{Th1}$
would become:
\begin{equation*}
\left[ \left( n^{-1}\tilde{T}_{[a_{n}t]}^{(n)},{c_{n}}^{-1}\tilde{H}
_{[a_{n}t]}^{(n)}\right) ,\,t\geq 0\right] \overset{(\mbox{\rm\tiny law})}{
\longrightarrow }(\tau ,H)\,,
\end{equation*}
which reduces to Theorem $3.2$ of $\cite{got}$ if the distribution of $
\tilde{S} ^{(n)}$ does not depend on $n.$
\end{rem}

\begin{proof}
We first recall Fristedt's formula for random walks, see \cite{do}, p.26. For every $\alpha>0$ and $\beta>0$,
\[1-\e\left(e^{-\alpha T^{(n)}_{1}-\beta H_{1}^{(n)}}\right)=\exp-\sum_{k=1}^{\infty}\frac{e^{-\alpha k}}k
\e\left(e^{-\beta S_k^{(n)}}:S_k^{(n)}>0\right)\,.\]
From this formula, we have
\begin{eqnarray*}
&&\e\left(e^{-\alpha n^{-1}T^{(n)}_{[a_n]}-\beta H_{[a_n]}^{(n)}}\right)=\e\left(e^{-\alpha n^{-1}T_1^{(n)}-\beta H^{(n)}_1}\right)^{[a_n]}\\
&=&\left(1-\exp-\sum_{k=1}^\infty\frac{1}ke^{-\alpha n^{-1}k}\e\left(e^{-\beta S_k^{(n)}}:S_k^{(n)}>0\right)\right)^{[a_n]}\\
&=&\left(1-\exp-\int_1^\infty\frac{1}{[s]}e^{-\alpha n^{-1}[s]}\e\left(e^{-\beta S_{[s]}^{(n)}}:S_{[s]}^{(n)}>0\right)\,ds\right)^{[a_n]}\\
&=&\left(1-\exp-\int_{1/n}^\infty\frac{n}{[nt]}e^{-\alpha n^{-1}[nt]}\e\left(e^{-\beta S_{[nt]}^{(n)}}:S_{[nt]}^{(n)}>0\right)\,dt\right)^{[a_n]}
\end{eqnarray*}
From the assumptions and Rogozin's criterion recalled in Remark \ref{rem1}, we have
\begin{eqnarray*}
\lim_{n\rightarrow+\infty}\int_{1/n}^\infty\frac{n}{[nt]}e^{-\alpha n^{-1}[nt]}\e\left(e^{-\beta S_{[nt]}^{(n)}}:S_{[nt]}^{(n)}>0\right)\,dt&=&\int_0^\infty\frac{e^{-\alpha t}}t\e\left(e^{-\beta X_t}:X_t>0\right)\,dt\\
&=&\infty\,,
\end{eqnarray*}
hence
\[-\ln\e\left(e^{-\alpha n^{-1}T^{(n)}_{[a_n]}-\beta H_{[a_n]}^{(n)}}\right)\sim
[a_n]\exp-\int_{1/n}^\infty\frac{n}{[nt]}e^{-\alpha n^{-1}[nt]}\e\left(e^{-\beta S_{[nt]}^{(n)}}:S_{[nt]}^{(n)}>0\right)\,dt\,.\]
From the expression of $a_n$ which is given in the statement of this theorem, the right and side of the above expression is
\begin{eqnarray*}
&&\exp\left(-\int_{1/n}^\infty\frac{n}{[nt]}e^{-\alpha n^{-1}[nt]}\e\left(e^{-\beta S_{[nt]}^{(n)}}:S_{[nt]}^{(n)}>0\right)\,dt
+\sum_{k=1}^\infty \frac1ne^{-k/n}\p(S_k^{(n)}>0)\right)\\
&=&\exp\int_{1/n}^\infty\frac{n}{[nt]}\e\left(e^{-n^{-1}[nt]}-e^{-\alpha n^{-1}[nt]-\beta S_{[nt]}^{(n)}}:S_{[nt]}^{(n)}>0\right)\,dt\,,
\end{eqnarray*}
which converges as $n$ goes to $+\infty$ toward
\[\exp\int_{0}^\infty\frac{1}{t}\e\left(e^{-t}-e^{-\alpha t-\beta X_t}:X_t>0\right)\,dt\\
=\kappa(\alpha,\beta)\,.\]
It is clear that the process $X$ drifts to $-\infty$ if and only if $S^{(n)}$ drifts to
$-\infty$ for all $n$ sufficiently large.  Suppose first that $X$ does not drift to $-\infty$.
The above convergence proves that the sequence $\left[ \left( n^{-1}T_{[a_{n}t]}^{(n)},H_{[a_{n}t]}^{(n)}\right) ,\,t\geq 0\right]$
converges in the sense of finite dimensional distributions toward $(\tau,H)$. We conclude that it converges weakly by applying Theorem 2.7
of Skorokhod \cite{sk}.

If $X$ drifts to $-\infty$, then the sequence $(T^{(n)},H^{(n)})$ and the process $(\tau,H)$ are obtained respectively from a
sequence of bivariate renewal processes, say $(\overline{T}^{(n)},\overline{H}^{(n)})$, and a bivariate subordinator, say $(\overline{\tau},\overline{H})$, all with infinite lifetime, by killing them respectively at
independent random times. It readily follows from the convergence of the characteristic exponents which is proved above that
$$\left[ \left( n^{-1}\overline{T}_{[a_{n}t]}^{(n)},\overline{H}_{[a_{n}t]}^{(n)}\right) ,\,t\geq 0
\right] \overset{(\mbox{\rm\tiny law})}{\longrightarrow }(\overline{\tau},\overline{H})$$ and that the independent
killing times of $\left( n^{-1}\overline{T}_{[a_{n}t]}^{(n)},\overline{H}_{[a_{n}t]}^{(n)}\right)$ converge in law to this of
$(\overline{\tau},\overline{H})$. As a straightforward consequence, the sequence of killed processes
$\left( n^{-1}T_{[a_{n}t]}^{(n)},H_{[a_{n}t]}^{(n)}\right)$ converges to  $(\tau,H)$, in the sense which is defined
in the preliminary section.
\end{proof}

\noindent Since $\tau$ is an increasing process, we derive from Theorem \ref
{Th1} and Theorem 7.2 of \cite{wh} that when $S^{(n)}$ converges in law to $X
$, the renormed process $\left(a_{n}^{-1}\Lambda _{\lbrack
nt]}^{(n)},\,t\ge0\right)$ converges in law to $(L_t,\,t\ge0)$. We actually establish
the following stronger result.

\begin{theorem}
\label{Theor1} Let $X$ be as in Theorem $\ref{Th1}$, and assume that
\begin{equation}
(S_{[nt]}^{(n)},\,t\geq 0)\overset{(
\mbox{\rm\tiny a.s.})}{\longrightarrow }(X_{t},\,t\geq 0).
\label{assume}
\end{equation}
Let $\Lambda ^{(n)}$ be the local time at its maximum of $S^{(n)}$. Then a
normed version of $\Lambda ^{(n)}$ converges uniformly in probability on
compacts sets towards $L$. More specifically, for all $t\geq 0$ and $
\varepsilon >0$,
\begin{equation}
\lim_{n\rightarrow \infty }\mathbb{P}\left( \sup_{s\in \lbrack
0,t]}|a_{n}^{-1}\Lambda _{\lbrack ns]}^{(n)}-L_{s}|>\varepsilon \right) =0\,,
\label{convuniprob}
\end{equation}
where $a_{n}$ is defined by expression $(\ref{a})$.
\end{theorem}

\noindent The proof of this theorem requires the two following lemmas. We
denote by $\pi^\tau$ and $\pi^H$ the L\'evy measures of $\tau$ and $H$.

\begin{lemma}
\label{lem0} The L\'evy measures $\pi^H$ has no atom whenever $X$ is not a
compound Poisson process. If moreover $0$ is regular for $(-\infty,0)$, then
the L\'evy measure $\pi^\tau$ has a density with respect to the Lebesgue
measure.
\end{lemma}

\begin{proof} Let us introduce some notations: we call $(\hat{\tau},\widehat{H})$  the ladder process
associated to $\widehat{X}=-X$ and we call $\widehat{U}$ the renewal measure of this bivariate subordinator.
The renewal measure of  $\widehat{H}$ is denoted by $U^{\widehat{H}}$ and the L\'evy measure of $X$
is denoted by $\Pi$.

From Vigon's \'equation amicale invers\'{e}e, see Vigon \cite{vi}, p.71, we have for all $x>0$
and $0<h<x$,
\begin{equation*}
{\pi ^{H}}[x-h,x)=\int_{0}^{\infty }U^{\widehat{H}}(dy){\Pi }[x+y-h,x+y)\,.
\end{equation*}
By monotone convergence, we get
\begin{equation*}
{\pi ^{H}}(\{x\})=\int_{0}^{\infty }U^{\widehat{H}}(dy){\Pi}(\{x+y\}),
\end{equation*}
and this is zero because there are countably many atoms of $\Pi,$ and
$U^{\widehat{H}}$ is diffuse when $X$ is not a compound Poisson process, see Proposition
1.15, Bertoin \cite{be}. This proves the first assertion.

Corollary 6, page 50 of \cite{do} asserts
that whenever $X$ is not a compound Poisson process, the L\'evy measure $\pi$ of $(\tau,H)$ is given by
\[\pi(dt,dh)=\int_{[0,\infty)}\widehat{U}(dt,dx)\,\Pi(dh+x)\,.\]
Then from Theorem 5 of \cite{ac}, under the additional
assumption that 0 is regular for $(-\infty,0)$,  we have for all $t>0$,
\[q_t(dx)\,dt=c\widehat{U}(dt,dx)\,,\]
where $c$ is a constant and $q_t(dx)$ is the entrance law of the measure of the excursions away from 0 of
the process $X$ reflected at its minimum. The second assertion is proved.
\end{proof}

\noindent The second lemma follows from Theorem \ref{Th1}, Lemma \ref{lem0}
and standard criterion on convergence of sums of independent of random
variables, see for instance \cite{fe}, so we omit its proof.

\begin{lemma}
\label{lem1} Define, for $0<a<b\leq \infty $, $0<c<\infty $ and $n\geq 1$,
\begin{eqnarray*}
\pi _{n}^{a,b} &=&\mathbb{P}(H_{1}^{(n)}\in (a,b])\,,\;m_{1}^{n,a}=\mathbb{E}
(H_{1}^{(n)}:H_{1}^{(n)}\leq a)\,,\;m_{2}^{n,a}=\mathbb{E}
((H_{1}^{(n)})^{2}:H_{1}^{(n)}\leq a) \\
\nu _{n}^{c} &=&\mathbb{P}(n^{-1}T_{1}^{(n)}>c)\,.
\end{eqnarray*}
Under the assumptions of Theorem $\ref{Th1}$, the following asymptotics
hold:
\begin{eqnarray*}
\lim_{n\rightarrow \infty }a_{n}\pi _{n}^{a,b} &=&\pi
^{H}(a,b]\,,\;\lim_{n\rightarrow \infty }a_{n}m_{1}^{n,a}=\delta
^{H}+\int_{0}^{a}x\,\pi ^{H}(dx)\,,\; \\
\text{and }\lim_{n\rightarrow \infty }a_{n}m_{2}^{n,a}
&=&\int_{0}^{a}x^{2}\,\pi ^{H}(dx)\,.
\end{eqnarray*}
If moreover $0$ is regular for $(-\infty ,0)$, then
\begin{equation*}
\lim_{n\rightarrow +\infty }a_{n}\nu _{n}^{c}=\pi ^{\tau }(c,\infty )\,.
\end{equation*}
\end{lemma}

\noindent
{\it Proof  of Theorem $\ref{Theor1}$}. We first observe that since
$(a_{n}^{-1}\Lambda _{\lbrack nt]}^{(n)},\,t\ge0)$ is a sequence of nondecreasing processes
which converges toward the continuous process $L$, in order to prove the uniform convergence in
$(\ref{convuniprob})$, it suffices to establish pointwise convergence in probability. However
this argument does not simplify the proof, so we deal directly with the uniform convergence.

We first treat the case where $\pi^H[0,\infty)<\infty$. Since we assumed that
0 is regular for $(0,\infty)$, we necessarily have $\delta_H>0$ and then
\begin{equation}\label{loctime}
\delta^H L_t= \lambda(M_s:s\le t)\,,\end{equation}
where $\lambda$ is the Lebesgue measure.
Let $M^{(n)}_k=\max_{0\le j\le k}S^{(n)}_j$, $k\ge0$ and for $a>0$, define the truncated past maxima of $S^{(n)}$
and $X$ respectively as
\[M^{n,a}_{[nt]}=M^{(n)}_{[nt]}-\sum_{s\in[0,t]}\Delta M^{(n)}_s\ind_{\{\Delta M^{(n)}_s>a\}}
\quad\mbox{and}\quad M^{a}_t= M_{t}-\sum_{s\in[0,t]}\Delta M_s\ind_{\{\Delta M_s>a\}}\,.\]
Since in this case, $M$ has only a finite number of jumps in each interval $[0,t]$, we have the
almost sure convergence
\begin{equation}\label{convas}
\lim_{n\rightarrow\infty} M^{n,a}_{[nt]}=M^{a}_t\,,\;\;\;\mbox{a.s.}\,.
\end{equation}
Moreover for the same reason and (\ref{loctime}), for all $a$ small enough we have
\begin{equation}\label{loctime2}
\delta^H L_t=M^a_t\,.\end{equation}
Then from (\ref{convas}) and (\ref{loctime2}), it is enough to prove that
\begin{equation}\label{mainlim}
\lim_{a\downarrow0}\limsup_{n\rightarrow\infty}\p\left(\sup_{t\in[0,1]}
|M^{n,a}_{[nt]}-\frac{\delta^H}{a_n}\Lambda^{(n)}_{[nt]}|>\varepsilon\right)=0\,.
\end{equation}
Note that for all $k$, $M^{n,a}(T^{(n)}_k)$ is the sum of $k$ i.i.d. random variables with mean $m_1^{n,a}$ and second
moment $m_2^{n,a}$ defined in Lemma \ref{lem1}. Hence for all $K>0$ and $\varepsilon>0$, from Kolmogorov's inequality,
\[\p\left(\max_{0\le j\le T^{(n)}_{Ka_n}}|M^{n,a}_j-m_1^{n,a}\Lambda_j^{(n)}|>\varepsilon  \right)\le
\frac{Ka_nm_2^{n,a}}{\varepsilon^2}\,.  \]
Now write the inequality
\begin{eqnarray*}
&&\p\left(\max_{0\le j\le T^{(n)}_{Ka_n}}|M^{n,a}_j-\frac{\delta^H}{a_n}\Lambda_j^{(n)}|>2\varepsilon  \right)\le\qquad\qquad\qquad\\
&&\qquad\qquad\qquad\qquad\p\left(\max_{0\le j\le T^{(n)}_{Ka_n}}|\frac{\delta^H}{a_n}\Lambda_j^{(n)}-m_1^{n,a}\Lambda_j^{(n)}|>2\varepsilon  \right)+
\frac{Ka_nm_2^{n,a}}{\varepsilon^2}\,.
\end{eqnarray*}
Then observe that the first term of the right hand side is nothing but $\ind_{\{|K\delta^H-m_1^{n,a}Ka_n|>\varepsilon\}}$
and from Lemma \ref{lem1}, $\lim_{a\rightarrow0}\lim_{n\rightarrow\infty}\ind_{\{|K\delta^H-m_1^{a,n}Ka_n|>\varepsilon\}}=0$.
From the same lemma, we have for the second term $\lim_{a\rightarrow0}\lim_{n\rightarrow\infty}a_nm_2^{n,a}=0$. Hence
\[\lim_{a\rightarrow0}
\lim_{n\rightarrow\infty}\p\left(\max_{0\le j\le T^{(n)}_{Ka_n}}|M^{n,a}_j-\frac{\delta^H}{a_n}\Lambda_j^{(n)}|>2\varepsilon  \right)=0\,.\]
Finally, write
\[\p\left(\sup_{t\in[0,1]}|M^{n,a}_{[nt]}-\frac{\delta^H}{a_n}\Lambda_{[nt]}^{(n)}|>\varepsilon  \right)\le
\p\left(\max_{0\le j\le T^{(n)}_{Ka_n}}|M^{n,a}_j-\frac{\delta^H}{a_n}\Lambda_j^{(n)}|>2\varepsilon  \right)+
\p(T^{(n)}_{Ka_n}<n)\,.\]
But from Lemma \ref{lem1}, we have
$\lim_{K\rightarrow+\infty}\lim_{n\rightarrow\infty}\p(T^{(n)}_{Ka_n}<n)=0$ and the conclusion follows in this case.

Now let us suppose that $\pi^H[0,\infty)=\infty$, and for $0<a<b<\infty$,
define the following approximations of the local times $L$ and $\Lambda^{(n)}$:
\[L^{a,b}_t=\#\{s\le t:\Delta M_s\in(a,b]\}\;\;\;\mbox{and}\;\;\;\Lambda^{n,a,b}_k=
\#\{j\le k:M^{(n)}_j+a<S^{(n)}_{j+1}\le M^{(n)}_j+b\}\,.\]
Since $L_t^{a,b}$ is a finite integer, it readily follows from the almost sure convergence of $S^{(n)}_{[n\cdot]}$ toward $X$ that
\begin{equation}\label{lim1}
\lim_{n\rightarrow+\infty}\Lambda^{n,a,b}_{[nt]}=L^{a,b}_t\,,\;\;\;\mbox{a.s.}\end{equation}
On the other hand, observe that $(L_{\tau_t}^{a,b},\,t\ge0)$ is a Poisson process with intensity $\pi^H(a,b]$.
Moreover from the hypothesis,
we have $\lim_{a\downarrow0}\pi^H(a,b]=+\infty$. So it follows from
the law of large numbers that for all $t>0$, $\lim_{a\rightarrow 0}\pi^H(a,b]^{-1}L_{\tau_t}^{a,b}=t$, a.s. From monotonicity, this
convergence can be strengthened to uniform convergence: for all $u>0$,
\[\lim_{a\rightarrow 0}\sup_{t\in[0,u]}|\pi^H(a,b]^{-1}L_{\tau_t}^{a,b}-t|=0\,,\;\;\;\mbox{a.s.}\]
Fix $\varepsilon>0$. For all $\eta>0$, we can  chose $u$ sufficiently large that $\p(\tau_u<1)<\eta/2$ and $a$
sufficiently small that $\p(\sup_{t\in[0,u]}|\pi^H(a,b]^{-1}L_{\tau_t}^{a,b}-t|>\varepsilon)<\eta/2$.
Then the inequality
\[\p(\sup_{t\in[0,1]}|\pi^H(a,b]^{-1}L_{t}^{a,b}-L_t|>\varepsilon)\le\p(\sup_{t\in[0,\tau_u]}|
\pi^H(a,b]^{-1}L_{t}^{a,b}-L_t|>\varepsilon)+\p(\tau_u<1)\]
allows us to obtain
\begin{equation}\label{lim2}
\lim_{a\rightarrow0}\p(\sup_{t\in[0,1]}|\pi^H(a,b]^{-1}L_{t}^{a,b}-L_t|>\varepsilon)=0\,.
\end{equation}
Note that for all $k$, $\Lambda^{n,a,b}(T^{(n)}_k)$ is the sum of $k$ independent Bernoulli random variables with mean $\pi_n^{a,b}$
defined in Lemma \ref{lem1}. Hence for all $K>0$, from Kolmogorov's inequality,
\[\p\left(\max_{0\le j\le T^{(n)}_{Ka_n}}\frac1{\pi^H(a,b]}|\Lambda^{a,b,n}_j-{\pi_n^{a,b}}\Lambda_j^{(n)}|>\varepsilon  \right)\le
\frac{Ka_n\pi_n^{a,b}}{\pi^H(a,b]^2\varepsilon^2}\,.  \]
Now write the inequality
\begin{eqnarray*}
&&\p\left(\max_{0\le j\le T^{(n)}_{Ka_n}}|\frac1{\pi^H(a,b]}\Lambda^{a,b,n}_j-\frac1{a_n}\Lambda_j^{(n)}|>\varepsilon  \right)\le\\
&&\p\left(\max_{0\le j\le T^{(n)}_{Ka_n}}|\frac{1}{a_n}\Lambda_j^{(n)}-\frac{\pi_n^{a,b}}{\pi^H(a,b]}\Lambda_j^{(n)}|>\varepsilon  \right)+
\frac{Ka_n\pi_n^{a,b}}{\pi^H(a,b]^2\varepsilon^2}\,.
\end{eqnarray*}
The first term of the right hand side is $\displaystyle\ind_{\{|K-Ka_n{\pi_n^{a,b}}/{\pi^H(a,b]}|>\varepsilon\}}$
and from Lemma \ref{lem1}, $\lim_{n\downarrow\infty}a_n{\pi_n^{a,b}}={\pi^H(a,b]}$, so this term converges to 0 for all $a$ and $b$.
The second term converges to $K/\varepsilon^2\pi^H(a,b]$ as $n$ tends to $\infty$. Since from the hypothesis we have
$\lim_{a\rightarrow0}\pi^H(a,b]=+\infty$, for all $b$, we conclude that
\[\lim_{a\rightarrow0}
\lim_{n\rightarrow\infty}\p\left(\max_{0\le j\le T^{(n)}_{Ka_n}}|
\frac1{\pi^H(a,b]}\Lambda^{a,b,n}_j-\frac1{a_n}\Lambda_j^{(n)}|>\varepsilon  \right)=0\,.\]
Finally, write
\begin{eqnarray*}&&\p\left(\sup_{t\in[0,1]}|
\frac1{\pi^H(a,b]}\Lambda^{a,b,n}_{[nt]}-\frac1{a_n}\Lambda_{[nt]}^{(n)}|>\varepsilon  \right)\le\\
&&\qquad\qquad\qquad\qquad\p\left(\max_{0\le j\le T^{(n)}_{Ka_n}}|
\frac1{\pi^H(a,b]}\Lambda^{a,b,n}_j-\frac1{a_n}\Lambda_j^{(n)}|>\varepsilon  \right)+
\p(T^{(n)}_{Ka_n}<n)\,.\end{eqnarray*}
The conclusion follows from (\ref{lim1}), (\ref{lim2}) and the fact that\\
$\lim_{K\rightarrow\infty}\lim_{n\rightarrow\infty}\p(T^{(n)}_{Ka_n}<n)=0$, which follows from Theorem \ref{Th1}.\hfill$\Box$\\
\noindent

\noindent When 0 is regular for $(-\infty,0)$, we may also define the local
time at the minimum of $X$ to be the local time at the maximum of $-X$. Let
us denote this process by $\widehat{L}$ and denote by $\widehat{\Lambda}
^{(n)}$ the local time at the maximum of the approximating sequence $-S^{(n)}
$. A straightforward consequence of the previous theorem is the following
result.

\begin{corollary}
\label{Coro3} Under the hypotheses of Theorem $\ref{Theor1}$,
\begin{equation*}
\left[\left(S^{(n)}_{[nt]},\frac1{a_n}\Lambda^{(n)}_{[nt]}\right),\,t\ge0
\right]\overset{\mbox{\rm\tiny (law)}}{\longrightarrow} [(X_t,L_t),\,t\ge0]
\,.
\end{equation*}
If in addition $0$ is regular for $(-\infty,0)$ then
\begin{equation*}
\left[\left(S^{(n)}_{[nt]},\frac1{a_n}\Lambda^{(n)}_{[nt]},\frac1{\hat{a}_n}
\widehat{\Lambda}^{(n)}_{[nt]}\right),\,t\ge0 \right]\overset{
\mbox{\rm\tiny
(law)}}{\longrightarrow} [(X_t,L_t,\widehat{L}_t),\,t\ge0] \,,
\end{equation*}
where $\hat{a}_{n}=\exp \left( \sum_{k=1}^{\infty }\frac{1}{k}e^{-k/n}
\mathbb{P}(S_{k}^{(n)}<0)\right)$.
\end{corollary}

\noindent

\begin{rem}
Let $S$ be a random walk whose law is in the domain of attraction of a
stable law. As an application of the previous corollary, in some instances,
we can compare the number of records of $S$ at its maximum with the number of
records at its minimum. More
precisely, if $a_n\sim \hat{a}_{n}$, then $\Lambda_{n}/\widehat{\Lambda}_{n}$
converges in law towards a non-degenerate random variable whereas if $\lim_n
\hat{a}_{n}/a_n=+\infty$, resp.~$0$, then $\Lambda_{n}/\widehat{\Lambda}_{n}$
goes to $+\infty$, resp.~$0$, in probability.
\end{rem}

\noindent When 0 is regular for $(-\infty,0)$, we denote by $(\hat{\tau},
\widehat{H})$ the strict ascending ladder process of $-X$ and for the
sequence of random walks $S^{(n)}$, we denote by $(\widehat{T}^{(n)},
\widehat{H}^{(n)})$ the strict ascending ladder height process of $-S$.
Another consequence of Theorem \ref{Theor1} is the following invariance
principle for both the ascending and descending ladder processes jointly.

\begin{corollary}
\label{Coro4} Let $X$ be any L\'{e}vy process such that $0$ is regular for
both the open half lines $(0,\infty )$ and $(-\infty,0)$ and assume that
some sequence of random walks $S^{(n)}$ converges in law toward $X$. Then
the process
\begin{equation*}
\left[\left( n^{-1}T_{[a_{n}t]}^{(n)},H_{[a_{n}t]}^{(n)},n^{-1}\widehat{T}_{[
\hat{a}_{n}t]}^{(n)}, \widehat{H}_{[\hat{a}_{n}t]}^{(n)}\right),\,t\ge0
\right]
\end{equation*}
converges toward the process
\begin{equation*}
[(\tau_t ,H_t,\hat{\tau}_t,\widehat{H}_t),\,t\ge0]
\end{equation*}
in the sense of finite dimensional distributions as $n\rightarrow\infty$,
where $a_n$ and $\hat{a}_n$ are defined respectively in Theorem $\ref{Th1}$
and Corollary $\ref{Coro3}$.
\end{corollary}

\begin{rem}
Note that in this case, we cannot conclude that weak convergence holds by
using Skorokhod's Theorem as in Theorem $\ref{Th1}$ since the quadrivariate
processes which are involved in Corollary $\ref{Coro4}$ do not have
independent increments.
\end{rem}

\begin{proof} By virtue of the Skorokhod representation theorem, there exists a sequence $\widetilde{S}^{(n)}$
(possibly defined on an enlarged probability space) such that
 for each $n$, $\widetilde{S}^{(n)}\ed S^{(n)}$ and such that $\widetilde{S}^{(n)}$ converges
almost surely toward $X$. Let $({\cal T}^{(n)},{\cal H}^{(n)})$ and $({\cal \widehat{T}}^{(n)},{\cal \widehat{H}}^{(n)})$
be respectively the strict ascending and the strict descending ladder processes of $\widetilde{S}^{(n)}$.

Recall that if a sequence of stochastic processes
converges almost surely on the Skorokhod space, then the sequence defined by the first passage time processes converges at
all continuity points of the limit process, see the remark after Theorem 7.1 in \cite{wh}. Moreover it is clear that the
 subordinators, $\tau$ and
$\hat{\tau}$ are a.s.~continuous at each $t\ge0$. Therefore, from Theorem \ref{Theor1} applied to $\widetilde{S}^{(n)}$ and
$-\widetilde{S}^{(n)}$, for fixed $t\ge0$, there exists a subsequence $(k_n)$ such that $k_n^{-1}{\cal T}_{[a_{k_n}t]}^{(k_n)}$ and
$k_n^{-1}\widehat{{\cal T}}_{[\hat{a}_{k_n}t]}^{(k_n)}$ converge almost surely toward $\tau_t$ and $\hat{\tau}_t$ respectively.

Since  $\tau_t$ and $\hat{\tau}_t$ are announceable stopping times in the filtration generated by $X$, it
follows from the quasi-left continuity of $X$ that this process is a.s.~continuous at times $\tau_t$ and $\hat{\tau}_t$, see
Ex.~3, Chap.~I in \cite{be}.

We deduce from the almost sure convergence of $\widetilde{S}^{(k_n)}$ toward $X$ that for every (possibly random) continuity point
$u$ of $X$, $\widetilde{S}^{(k_n)}_{[k_nu]}$ converges almost surely to $X_u$, see \cite{bi}, p.112. Therefore the sequence
\begin{eqnarray*}
&&\left(k_n^{-1}{\cal T}^{(k_n)}_{[{a_{k_n}}t]},
{\cal H}^{(k_n)}_{[{a}_{k_n}t]},k_n^{-1}{\cal \widehat{T}}^{(k_n)}_{[\hat{a}_{k_n}t]},
{\cal \widehat{H}}^{(k_n)}_{[\hat{a}_{k_n}t]}\right)\\
&&\qquad\qquad\qquad\qquad=\left(k_n^{-1}{\cal T}^{(k_n)}_{[{a_{k_n}}t]},\widetilde{S}^{(k_n)}\left({\cal T}^{(k_n)}_{[{a_{k_n}}t]}\right),k_n^{-1}{\cal \widehat{T}}^{(k_n)}_{[\hat{a}_{k_n}t]},-\widetilde{S}^{(k_n)}\left(\widehat{{\cal T}}^{(k_n)}_{[\hat{a}_{k_n}t]}\right)\right)
\end{eqnarray*}
converges almost surely toward $(\tau_t,X(\tau_t),\hat{\tau}_t,-X(\hat{\tau}_t))=(\tau_t,H_t,\hat{\tau}_t,\widehat{H}_t)$,
as $n\rightarrow\infty$. This almost sure convergence
is easily extended to the multidimensional case, i.e. there is a subsequence $(k'_n)$ such that
it holds simultaneously at any sequence of times $t_1,\dots,t_j$. So we have proved that the variables
$\left\{ n^{-1}{\cal T}_{[a_{n}t_i]}^{(n)},{\cal H}_{[a_{n}t_i]}^{(n)},n^{-1}\widehat{{\cal T}}_{[\hat{a}_{n}t_i]}^{(n)},
\widehat{{\cal H}}_{[\hat{a}_{n}t_i]}^{(n)},\,i=1,\dots j\right\}$ converge in probability,
and  we conclude from the identity in law
\[({\cal T}^{(n)},{\cal H}^{(n)},{\cal \widehat{T}}^{(n)},{\cal \widehat{H}}^{(n)})\ed
(T^{(n)},H^{(n)},\widehat{T}^{(n)},\widehat{H}^{(n)})\,,\]
which holds for each $n$, as a consequence of the identity $\widetilde{S}^{(n)}\ed S^{(n)}$.\end{proof}

\noindent Now we suppose that there is a sequence of L\'{e}vy processes $
X^{(n)}$, all of which satisfy the same hypothesis as $X$, i.e. 0 is regular
for $(0,\infty)$. Call $L^{(n)}$ the version of the local time of $X^{(n)}$
at its maximum, as it is defined in Section \ref{prelim}.

\begin{theorem}
\label{th2} Suppose that as $n$ tends to $\infty $,
\begin{equation*}
X^{(n)}\overset{\mbox{\rm\tiny a.s.}}{\longrightarrow }X\,.
\end{equation*}
Then the sequence of local times $L^{(n)}$ converges uniformly on compact
sets in probability toward $L$, i.e. for all $t>0$ and $\varepsilon >0$,
\begin{equation*}
\lim_{n\rightarrow \infty }\mathbb{P}\left( \sup_{s\in \lbrack
0,t]}|L_{t}^{(n)}-L_{t}|>\varepsilon \right) =0\,.
\end{equation*}
\end{theorem}

\begin{proof} For each $n$, we define a sequence of random walks $(S^{n,k},\,k\ge0)$ whose paths are embedded
in those of $X^{(n)}$ as follows:
\[S^{n,k}_j=X^{(n)}_{j/k}\,,\;j\ge0\,.\]
Then we may readily check that for each $n$, as $k$ tends to $\infty$,
\[(S^{n,k}_{[kt]},\,t\ge0)\stackrel{\mbox{\rm\tiny a.s.}}{\longrightarrow}X^{(n)}\,.\]
Call $\Lambda^{n,k}$ the local time at its maximum of $S^{n,k}$ as it is defined for $S^{(n)}$ in (\ref{loctimerw}).
From Theorem \ref{Theor1}, we have for all $n\ge1$, $t\ge0$ and $\varepsilon>0$,
\[\lim_{k\rightarrow\infty}\p\left(\sup_{s\in[0,t]}|\frac1{a_k^n}\Lambda^{n,k}_{[ns]}-L_s^{(n)}|>\varepsilon\right)=0\,,\]
where $\ln a_k^n=\sum_{j=1}^\infty \frac1je^{-j/k}\p(S_j^{n,k}>0)$. We can chose a sequence of integers $(k_n)_{n\ge1}$ such that
\[\lim_{n\rightarrow\infty}\p\left(\sup_{s\in[0,t]}|\frac1{a_{k_n}^n}\Lambda^{n,{k_n}}_{[ns]}-L_s^{(n)}|>\varepsilon\right)=0\]
and, as $n$ goes to $\infty$,
\[ (S_{[k_nt]}^{n,k_n},\,t\ge0)\stackrel{\mbox{\rm\tiny a.s.}}{\longrightarrow}X\,.\]
Hence by applying again Theorem \ref{Theor1}, we have
\[\lim_{n\rightarrow\infty}\p\left(\sup_{s\in[0,t]}|\frac1{a_{k_n}^n}\Lambda^{n,k_n}_{[ns]}-L_s|>\varepsilon\right)=0\,,\]
which allows us to conclude.
\end{proof}
\noindent It is clear that the equivalent results to Corollaries \ref{Coro3}
and \ref{Coro4} are also valid in the setting of Theorem \ref{th2}, that is
replacing the approximating sequence $S^{(n)}$ by the sequence $X^{(n)}$.

\section{Applications to conditioned processes}

In this section we will prove that if a sequence $S^{(n)}$ of random walks
converges weakly toward a L\'evy process $X$, then the sequence $S^{(n)}$
conditioned to stay positive also converges weakly toward $X$ conditioned to
stay positive. For simplicity in the statements and proofs, we will always
suppose that $S^{(n)}$ and $X$ do not drift to $-\infty$ and that for $X$,
the state $0$ is regular for both $(-\infty,0)$ and $(0,\infty)$.

We first define $S^{(n)}$ and $X$ conditioned to stay positive on the whole
time interval $[0,\infty)$. Let $V^{(n)}(x)=\sum_{k\ge0}\mathbb{P}(\widehat{H
}^{(n)}_k\le x)$, $x\ge0$ be the renewal function of $\widehat{H}^{(n)}$ and
let $S^{(n)*}$ be the process $S^{(n)}$ killed when it enters the negative
half-line. The Markovian transition function which is given by
\begin{equation*}
q^\uparrow(x,dy)=\frac{V^{(n)}(y)}{V^{(n)}(x)}\mathbb{P}(S_{k+1}^{(n)*}\in
dy\,|S_k^{(n)*}=x)\,,
\end{equation*}
for $x>0$, $y>0$ if $k\ge1$, and for $x\ge0$, $y>0$ if $k=0$, characterizes
the law of an $h$-process of $S^{(n)*}$ which is called \textit{the law of $
S^{(n)}$ conditioned to stay positive.} Similarly, denote by $X^*$ the
L\'evy process $X$ killed when it enters the negative halfline. Suppose that
0 is regular for $(-\infty,0)$ and let $h(x)$ be the renewal function of the
subordinator $\widehat{H}$, i.e. $h(x)=\mathbb{E}\left(\int_0^\infty
\mbox{\rm 1\hspace{-0.04in}I}_{\{\widehat{H}_t\le x\}}\,dt\right)$, then the
Markovian semigroup
\begin{equation*}
p^\uparrow_t(x,dy)=\frac{h(y)}{h(x)}\mathbb{P}(X^*_{t+s}\in
dy\,|X_s^*=x)\,,\quad\mbox{for $x,y>0$ and $s,t>0$}
\end{equation*}
is that of the L\'evy process $X$ \textit{conditioned to stay positive.} For
$x=0$, this semigroup admits a unique entrance law which is specified in
terms of the measure of the excursions above the minimum of the process $X$.
We refer to \cite{bd}, \cite{ch} and \cite{cd} for a more complete account
on random walks and L\'evy processes conditioned to stay positive.

The proof of our invariance principle bears upon a pathwise construction of $
S^{(n)}$ and $X$ conditioned to stay positive which is due to Tanaka and
Doney, see \cite{do1} and \cite{do}, p.91. Let us briefly recall it both in
discrete time and in continuous time. For $k\ge0$, call $e^{(k)}$ the $k$-th
excursion of the reflected process $M^{(n)}-S^{(n)}$:
\begin{equation*}
e_i^{(k)}=(M^{(n)}-S^{(n)})_{T_k^{(n)}+i}\,,\;\;\;\;0\leq i\le
T_{k+1}^{(n)}-T_k^{(n)}\,,
\end{equation*}
and denote by $\hat{e}^{(k)}$ the time reversal of $e^{(k)}$, i.e.
\begin{equation*}
\hat{e}_i^{(k)}=H^{(n)}_{k+1}-S^{(n)}(T^{(n)}_{k+1}-i)\,,\;\;\;\;0\leq i\le
T_{k+1}^{(n)}-T_k^{(n)}\,.
\end{equation*}
The process $S^{(n)\uparrow}$ which is obtained from the concatenation of $
\hat{e}^{(0)},\hat{e}^{(1)},...$, i.e:
\begin{equation}  \label{td}
S^{(n)\uparrow}_i=H^{(n)}_k+\hat{e}^{(k)}_{i-T_k^{(n)}}, \quad\mbox{if}\quad
T_k^{(n)}\le i\le T_{k+1}^{(n)}\,,
\end{equation}
has the law of $S^{(n)}$ conditioned to stay positive. A similar
construction in continuous time has been obtained in \cite{do}: for $t>0$,
let
\begin{equation*}
g(t)=\sup\{s<t:X_s=M_s\}\quad\mbox{and}\quad d(t)=\inf\{s>t:X_s=M_s\}\,,
\end{equation*}
then the process defined by
\begin{equation}  \label{tdc}
X_0^\uparrow=0\quad\mbox{and}\quad X_t^\uparrow=M_{d(t)}+\mbox{\rm 1\hspace{-0.04in}I}
_{\{d(t)>g(t)\}}(M-X)_{(d(t)+g(t)-t)-}\,,\quad t>0
\end{equation}
has the law of $X$ conditioned to stay positive.

Let us also define the local time at the future minimum of $S^{(n)\uparrow}$
and $X^\uparrow$. The first of these processes is simply the counting
process defined by $\underline{\underline{\Lambda}}^{(n)}_0=0$ and for $k\ge1$,
\begin{equation*}
\underline{\underline{\Lambda}}^{(n)}_k= \#\{j\in \{1,\dots
,k\}:S_{j}^{(n)\uparrow}=\min_{i\geq
j}S_{i}^{(n)\uparrow},\,S_j^{(n)\uparrow}<S_{j+1}^{(n)\uparrow}\}\,.
\end{equation*}
Recall that in continuous time the set $\{t:X_t^\uparrow=\inf_{s\ge
t}X^\uparrow_s\}$ is regenerative so that we may define on this set a local
time $\underline{\underline{L}}$, see \cite{ch}, p.44. This local time is
unique up to a normalizing constant and we will normalize it by $\mathbb{E}
\left(\int_0^\infty e^{-t}\,d\underline{\underline{L}}_t\right)=1$. One
easily derives from the above pathwise constructions the identities
\begin{equation*}
\{j\ge1:S_{j}^{(n)\uparrow}=\min_{i\geq
j}S_{i}^{(n)\uparrow},\,S_j^{(n)\uparrow}<S_{j+1}^{(n)\uparrow}\}
=\{j\ge1:S_{j-1}<S_j,\,S_{j}=\max_{i\leq j}S_{i}\}
\end{equation*}
and $\{t:X_t^\uparrow=\inf_{s\ge t}X^\uparrow_s\}=\{t:X_t=\sup_{s\le t}X_s\}$.
In particular, we have
\begin{equation}  \label{idloc}
\underline{\underline{\Lambda}}^{(n)}=\Lambda^{(n)}\quad\mbox{and}\quad
\underline{\underline{L}}=L\,,\quad\mbox{a.s.}
\end{equation}

The following theorem has been partially obtained in the particular setting
of stable processes in \cite{cc}, see Theorem 1.1.

\begin{theorem}
\label{th3} Suppose that some sequence of random walks $S^{(n)}$ converges
almost surely toward $X$. Recall the definition of $a_n$ from Theorem $\ref{Th1}$.

\begin{itemize}
\item[$1.$] The sequence of processes $(S^{(n)\uparrow}_{[nt]},\,t\ge0)$
converges almost surely toward $X^\uparrow$.

\item[$2.$] The sequence $[(S^{(n)\uparrow}_{[nt]},a_n^{-1}\underline{
\underline{\Lambda}}^{(n)}_{[nt]}),\,t\ge0]$ converges in probability toward
$(X^\uparrow,\underline{\underline{L}})$.
\end{itemize}

Consequently, if some sequence of random walks $S^{(n)}$ converges weakly
toward $X$, then the sequence $[(S^{(n)\uparrow}_{[nt]},a_n^{-1}\underline{
\underline{\Lambda}}^{(n)}_{[nt]}),\,t\ge0]$ converges weakly toward $
(X^\uparrow,\underline{\underline{L}})$.
\end{theorem}

\noindent In the second part of Theorem \ref{th3}, convergence in
probability means that each coordinate converges in probability with respect
to some distance which generates Skorokhod topology on the space $\mathcal{D}
([0,\infty))$. But more particularly, from part 1., the first coordinate
converges almost surely, whereas the second coordinate converges uniformly
in probability on compact sets, in the sense which has been defined in
Theorem \ref{Theor1}. The result displayed in Theorem \ref{th3} holds in the
very general case, although for simplicity in its statement and proof we
restrict ourself to the case where 0 is regular.

The time reversal relationships between $X$ and $X^\uparrow$ and between $
S^{(n)}$ and $S^{(n)\uparrow}$ which are presented below, in Theorem \ref
{th4} and Lemma \ref{timereversal},  are required for the proof of Theorem
\ref{th3}. Let us denote by $U^{(n)}$ and $\sigma$ respectively, the (right
continuous) inverses of $\underline{\underline{\Lambda}}^{(n)}$ and $
\underline{\underline{L}}$, i.e:
\begin{equation*}
U^{(n)}_k=\min\{i:\underline{\underline{\Lambda}}^{(n)}_i=k\}\,,\;\;k\ge0
\quad\mbox{and}\quad\sigma_t=\inf\{s:\underline{\underline{L}}
_s>t\}\,,\;\;t\ge0\,.
\end{equation*}
We also set
\begin{equation*}
G_k^{(n),\uparrow}=\max\{j\le k:S^{(n)\uparrow}_j= \inf_{i\ge
j}S_i^{(n),\uparrow}\}\quad\mbox{and}\quad g^\uparrow_t=\sup\{s\le
t:X_s^\uparrow=\inf_{u\ge s}X^\uparrow_u\} \,,
\end{equation*}
and $G_k^{(n)}=\max\{j\le k:M_j^{(n)}=S_j^{(n)}\}$.

\begin{theorem}
\label{th4} The following time reversal relationships hold between $X$ and
$X^\uparrow$$:$

\begin{itemize}
\item[$1.$] For all $t>0$, the law of the process $[(X_{\tau_t}-X_{(
\tau_t-s)-},L_{\tau_t}-L_{\tau_t-s}),\,0\le s<\tau_t]$ is the same as that
of the process $[(X^\uparrow_{s},\underline{\underline{L}}_s),\,0\le
s<\sigma_t]$.

\item[$2.$] For all $t>0$, the law of the process $
[(X_{g(t)}-X_{(g(t)-s)-},L_{g(t)}-L_{g(t)-s}),\,0\le s\le g(t)]$ $($with the
convention that $0-=0$$)$ is the same as that of the process $
[(X^\uparrow_{s},\underline{\underline{L}}_t),\,0\le s\le g^\uparrow_t]$.
\end{itemize}
\end{theorem}

\noindent Note that in the above statement, we have $X_{g(t)}=M_t$ and $
L_{\tau_t}=t$, almost surely. Part~1.~of this theorem is Lemma 4.3 of
Duquesne \cite{du}. The case where these processes have no positive jumps,
is treated in Theorem VII.18 of \cite{be}. It generalizes a well known
transformation between Brownian motion and the three dimensional Bessel
process due to Williams. Here we show that this result can easily be derived
from simple arguments involving Tanaka-Doney's transformation. Our next
lemma states the discrete time counterpart of Theorem \ref{th4}. Its proof
is very similar to that of Theorem \ref{th4}, hence we will only prove the
discrete time case.

\begin{lemma}
\label{timereversal} For any $k\geq 1$,

\begin{itemize}
\item[$1.$] the law of the process $
[(S^{(n)}_{T_k^{(n)}}-S^{(n)}_{T_k^{(n)}-i},k-\Lambda^{(n)}(T_k^{(n)}-i)),
\,0\le i\le T_k^{(n)}]$ is the same as that of the process $
[(S^{(n)\uparrow}_i,\underline{\underline{\Lambda}}^{(n)}_i),\,0\le i\le
U^{(n)}_k)]$,

\item[$2.$] the law of the process $
[(S^{(n)}_{G_k^{(n)}}-S^{(n)}_{G_k^{(n)}-i},\Lambda^{(n)}_{G_k^{(n)}}-
\Lambda^{(n)}_{G_k^{(n)}-i}),\,0\le i\le G_k^{(n)}]$ is the same as that of
the process $[(S^{(n)\uparrow}_i,\underline{\underline{\Lambda}}
^{(n)}_i),\,0\le i\le G^{(n)\uparrow}_k]$.
\end{itemize}
\end{lemma}

\begin{proof}  From the transformation which is recalled in
(\ref{td}), the process $S^{(n)\uparrow}$ is the concatenation of
the time reversed excursions $\hat{e}^{(0)},\hat{e}^{(1)},\dots$. It is clear that the times where this process reaches its
future minimum occur at the end of each of these reversed excursions.  Therefore $T_k^{(n)}=U^{(n)}_k$, a.s.~and
the concatenation of the $k$ excursions
$\hat{e}^{(0)},\hat{e}^{(1)},\dots,\hat{e}^{(k)}$ is the process $(S^{(n)\uparrow}_i,\,0\le i\le U^{(n)}_k)$.

From the Markov property, these excursions
are i.i.d., so that the concatenation of $\hat{e}^{(0)},\hat{e}^{(1)},\dots,\hat{e}^{(k)}$ has the same law as the
concatenation of $\hat{e}^{(k)},\hat{e}^{(k-1)},\dots,\hat{e}^{(1)}$. But the latter concatenation is precisely
the process $(S^{(n)}_{T_k^{(n)}}-S^{(n)}_{T_k^{(n)}-i},\,0\le i\le T_k^{(n)})$. The same reasoning justifies the identity on the
second coordinate.

The second part of the statement follows from the same arguments together with the identity $G_k^{(n)}=G_k^{(n)\uparrow}$
which holds for all $k\ge0$.
\end{proof}

\noindent Actually in the proof of Theorem \ref{th3} we will only use the
second part of Theorem \ref{th4} which says that the returned pre-maximum
part of $X$ before time $t$ has the same law as $X^\uparrow$ up to its last
passage time at the future minimum before $t$. However, in order to avoid
the need to justify an invariance principle for returned processes, we will
reformulate this identity in law in terms of the post-minimum process.
\newline

\noindent {\it Proof of Theorem $\ref{th3}$}. From identity (\ref{idloc})
and Theorem \ref{Theor1}, we only have to prove part~1.~of Theorem \ref{th3}.
Define
\begin{equation*}
K^{(n)}_j=\sup\{i\le j:S^{(n)}_i=\min_{l\le i} S^{(n)}_l\}\quad\mbox{and}
\quad k(t)=\sup\{s\le t:X_s=\inf_{u\le s}X_u\}\,.
\end{equation*}
From time reversal properties of $S^{(n)}$ and $X$, we have:
\begin{equation*}
\left(S^{(n)}_{K^{(n)}_k+i}-S^{(n)}_{K^{(n)}_k},\,0\le i\le
k-K^{(n)}_{k}\right)\overset{(d)}{=}
\left(S^{(n)}_{G_k^{(n)}}-S^{(n)}_{G_k^{(n)}-i},\,0\le i\le G_k^{(n)}\right)
\end{equation*}
and
\begin{equation*}
(X_{k(t)+s}-X_{k(t)},\,0\le s\le t- k(t))\overset{(d)}{=}
(X_{g(t)}-X_{(g(t)-s)-},\,0\le s\le g(t))\,.
\end{equation*}
(Recall the convention: $0-=0$). Since 0 is regular for both $(-\infty,0)$
and $(0,\infty)$, the time $k(t)$ is a continuity point of $X$, hence the
almost sure convergence of $S^{(n)}$ toward $X$ implies that $\lim_n
n^{-1}K^{(n)}_{[nt]}=k(t)$, a.s.~for all $t\ge0$. Then recall from the
preliminary section our definition of the a.s. convergence of stochastic
processes with finite lifetime. We clearly have the almost sure convergence
of the sequence of processes
\begin{equation*}
Y^{(n)}=\left(S^{(n)}_{K^{(n)}_{[nt]}+[ns]}-S^{(n)}_{K^{(n)}_{[nt]}},\,0\le
s\le n^{-1}([nt]-K^{(n)}_{[nt]})\right)
\end{equation*}
toward the process $(X_{k(t)+s},\,0\le s\le t-k(t))$. From Lemma \ref
{timereversal} and the time reversal property of $S^{(n)}$, the sequence
$Y^{(n)}$, $n\ge0$ has the same law as the sequence
\begin{equation*}
\left(S^{(n)\uparrow}_{[ns]},\,0\le s\le
n^{-1}G^{(n)\uparrow}_{[nt]}\right)\,.
\end{equation*}
Therefore the sequence $\left(S^{(n)\uparrow}_{[ns]},\,0\le s\le
n^{-1}G^{(n)\uparrow}_{[nt]}\right)$ converges almost surely toward the
process $(X^\uparrow_{s},\,0\le s\le g^\uparrow_t)$.\newline

Let $(t_k)$ be an increasing sequence of positive reals which tends to $
\infty$.  We deduce from the
above convergence that for each $k$, $\lim_{n\rightarrow\infty}n^{-1}G^{(n)
\uparrow}_{[nt_k]}=g^\uparrow(t_k)$, a.s. and more generally,
\begin{equation*}
\left(S^{(n)\uparrow}_{[ns]}\mbox{\rm
1\hspace{-0.04in}I}_{\left\{n^{-1}G^{(n)\uparrow}_{[nt_k]}\le 1<n^{-1}
G^{(n)\uparrow}_{[nt_{k+1}]}\right\}},\,0\le s\le 1\right)
\end{equation*}
converges a.s. toward $(X^\uparrow_s\mbox{\rm
1\hspace{-0.04in}I}_{\{g^\uparrow(t_k)\le 1<g^\uparrow(t_{k+1})\}},\,0\le s\le 1)$.
Since all processes $S^{(n)\uparrow}$ and $X^\uparrow$ drift to $+\infty$,
we have $\lim_{k\rightarrow\infty}G^{(n)\uparrow}_{[nt_k]}=\infty$ and $
\lim_{k\rightarrow\infty}g^\uparrow(t_k)=+\infty$, a.s., so that with $t_0=0$, we have
\begin{eqnarray*}
\left(S^{(n)\uparrow}_{[ns]},\,0\le s\le1\right)&=&\left(\sum_{k\ge0}S^{(n)\uparrow}_{[ns]}\mbox{\rm
1\hspace{-0.04in}I}_{\left\{n^{-1}G^{(n)\uparrow}_{[nt_k]}\le 1<n^{-1}
G^{(n)\uparrow}_{[nt_{k+1}]}\right\}},\,0\le s\le 1\right)\quad\mbox{and}\\
(X^\uparrow_s,\,0\le s\le 1)&=&\left(\sum_{k\ge0}X^\uparrow_s\mbox{\rm
1\hspace{-0.04in}I}_{\{g^\uparrow(t_k)\le 1<g^\uparrow(t_{k+1})\}},\,0\le s\le 1\right)\,.
\end{eqnarray*}
But almost surely there is $k$ and $n_0$ such that for all $n\ge n_0$,  the processes
on the right hand sides of the two equalities above are respectively equal to
$$\left(S^{(n)\uparrow}_{[ns]}\mbox{\rm
1\hspace{-0.04in}I}_{\left\{n^{-1}G^{(n)\uparrow}_{[nt_k]}\le 1<n^{-1}
G^{(n)\uparrow}_{[nt_{k+1}]}\right\}},\,0\le s\le 1\right)$$ and $(X^\uparrow_s\mbox{\rm
1\hspace{-0.04in}I}_{\{g^\uparrow(t_k)\le 1<g^\uparrow(t_{k+1})\}},\,0\le s\le 1)$.
Therefore $\left(S^{(n)\uparrow}_{[ns]},\,s\ge0\right)$ converges
toward $(X^\uparrow_s,\,s\ge0)$ on the space  $\mathcal{D}([0,1])$. The same arguments
holds on each space $\mathcal{D}([0,t])$, $t>0$ so we deduce the convergence
on $\mathcal{D}([0,\infty))$ from Theorem 16.7 in \cite{bi} as recalled in
the preliminary section.\hfill$\Box$\\

We now define $S^{(n)}$ and $X$ conditioned to stay positive respectively on
$\{0,1,\dots,k\}$ and $[0,t]$, where $k$ and $t$ are deterministic. Let $
\mathcal{C}_k^{(n)}=\{S^{(n)}_1\ge0,\dots,S^{(n)}_k\ge0\}$ then we denote by
$S^{(n,k)}$ a process whose law is defined on $\{0,1,\dots,k\}$ by $
S_0^{(n,k)}=0$ and
\begin{equation*}
\mathbb{P}(S^{(n,k)}_1\in dx_1,\dots,S^{(n,k)}_k\in dx_k)= \mathbb{P}
(S^{(n)}_1\in dx_1,\dots,S^{(n)}_k\in dx_k\,|\mathcal{C}_k^{(n)})\,.
\end{equation*}
It clearly follows from the definitions that this law is absolutely
continuous with respect to the law of $S^{(n)\uparrow}$: for $
x_1>0,\dots,x_k>0$,
\begin{equation}  \label{rel1}
\mathbb{P}(S^{(n,k)}_1\in dx_1,\dots,S^{(n,k)}_k\in dx_k)= \frac{1}{\mathbb{P
}(\mathcal{C}_k^{(n)})V^{(n)}(x_k)}\mathbb{P}(S^{(n)\uparrow}_1\in
dx_1,\dots,S^{(n)\uparrow}_k\in dx_k)\,.
\end{equation}
See also (3.2) in \cite{cc}. The process $S^{(n,k)}$ is called the (discrete
time) \textit{meander} with length~$k$.

The definition of the analogous conditional law in continuous time requires
some care since the set $\{X_t\ge0:t\in[0,1]\}$ has always probability 0
when 0 is regular for $(-\infty,0)$.

\begin{lemma} 
For $x_1>0,\dots,x_j>0$ and $t_1,\dots,t_j\in[0,1]$, we have
\begin{eqnarray*}
&&\lim_{x\rightarrow0}\mathbb{P}_x(X_{t_1}\in dx_1,\dots,X_{t_j}\in
dx_j\,|\,X_t>0,t\in[0,1]) \\
&&\qquad\qquad\qquad\qquad\qquad\qquad =\frac{1}{\beta h(x_j)}\mathbb{P}
(X^\uparrow_{t_1}\in dx_1,\dots,X^\uparrow_{t_j}\in dx_j)\,,
\end{eqnarray*}
where $\beta=\mathbb{E}(h(X_1^\uparrow)^{-1})$.
\end{lemma}

\begin{proof} This is a direct application of Corollary 1 in \cite{cd}, see also \cite{cd1}.
\end{proof}
\noindent Clearly the weak limit obtained in this lemma defines a unique
probability measure on the space $\mathcal{D}([0,1])$. We will denote by $X^+
$ a process with this law, i.e. for $x_1>0,\dots,x_j>0$ and $
t_1,\dots,t_j\in[0,1]$,
\begin{equation}  \label{rel2}
\mathbb{P}(X^+_{t_1}\in dx_1,\dots,X^+_{t_j}\in x_j)= \frac{1}{\beta h(x_j)}
\mathbb{P}(X^\uparrow_{t_1}\in dx_1,\dots,X^\uparrow_{t_j}\in dx_j)\,.
\end{equation}
This process is called the \textit{meander} with length 1.

\begin{lemma}
\label{harmfunct} Assume that $S^{(n)}$ converges weakly to $X$. Recall the definition of the renewal function $
V^{(n)}(x)=\sum_{k\ge0}\mathbb{P}(\widehat{H}^{(n)}_k\le x)$, for $x\ge0$.

\begin{itemize}
\item[$1.$] Let $\pi^{\hat{\tau}}$ be the L\'evy measure of the ladder time
process $\hat{\tau}$, then
\begin{equation*}
\lim_{n\rightarrow+\infty}\hat{a}_n\mathbb{P}(\mathcal{C}_n^{(n)})= \pi^{
\hat{\tau}}(1,\infty)\,.
\end{equation*}

\item[$2.$] The sequence of functions $\mathbb{P}(\mathcal{C}
_n^{(n)})V^{(n)}(x)$ converges uniformly on compacts sets toward $\gamma h(x)=
\gamma\mathbb{E}\left(\int_0^\infty\mbox{\rm 1\hspace{-0.04in}I}_{\{\widehat{H}
_t\le x\}}\,dt\right)$, with $\gamma=\pi^{\hat{\tau}}(1,\infty)$.
\end{itemize}
\end{lemma}

\begin{proof} To prove the first part, it suffices to note that
$\p({\cal C}_n^{(n)})=\p(n^{-1}\widehat{T}_1^{(n)}>1)$ and to apply
Lemma \ref{lem1}. To prove the second part, observe that from the hypothesis,
Theorem \ref{Th1} and dominated convergence, we have for every $x\ge0$,
\[\lim_{n\rightarrow\infty}\int_0^\infty\p(\widehat{H}^{(n)}_{[\hat{a}_nt]}\le x)\,dt=\lim_{n\rightarrow\infty}
\hat{a}_n^{-1}V^{(n)}(x)=h(x)\,.\]
Then the result follows from part 1., the fact that $V^{(n)}(x)$ is a sequence of increasing functions and the
continuity of $h$.
\end{proof}

\noindent The following invariance principle for the meander has been obtained in the case where all
$S^{(n)}$ have the same law (in particular $X$ is stable) in \cite{bo} and \cite{do0}.

\begin{theorem}
\label{th6} Suppose that some sequence of random walks $S^{(n)}$ converges
weakly toward $X$. The sequence of discrete meanders $(S^{(n,n)}_{[nt]},\,0
\le t\le1)$ converges weakly toward the meander $X^+$.
\end{theorem}

\begin{proof} We will prove that for all continuous and bounded functionals $F$ on ${\cal D}([0,1])$,
\[\e\left(F(S^{(n,n)}_{[nt]},\,0\le t\le1)\right){\longrightarrow}
\e\left(F(X^+_t,\,0\le t\le1)\right),
\qquad\mbox{as $n\rightarrow\infty$.}\]
From the absolute continuity
relations (\ref{rel1}) and (\ref{rel2}), it suffices to prove that
\begin{eqnarray*}
&&\e\left(\frac1{\p({\cal C}_n^{(n)})V^{(n)}(S_n^{(n)\uparrow})}F(S^{(n)\uparrow}_{[nt]},\,0\le t\le1)\right)\\
&&\qquad\qquad\qquad\qquad\qquad\qquad{\longrightarrow}\;\;
\e\left(\frac1{\beta h(X_1^\uparrow)}F(X^\uparrow_t,\,0\le t\le1)\right),
\qquad\mbox{as $n\rightarrow\infty$.}
\end{eqnarray*}
For $\eta>0$, write
\begin{eqnarray*}
&&\left|\e\left(\frac1{\p({\cal C}_n^{(n)})V^{(n)}(S_n^{(n)\uparrow})}F(S^{(n)\uparrow}_{[nt]},\,0\le t\le1)\right)-\e\left(\frac1{\gamma h(X_1^\uparrow)}F(X^\uparrow_t,\,0\le t\le1)\right)\right|\\
&&\le\left|\e\left(\frac1{\p({\cal C}_n^{(n)})V^{(n)}(S_n^{(n)\uparrow})}
\ind_{\{S^{(n)\uparrow}_n\ge\eta\}}F(S^{(n)\uparrow}_{[nt]},\,0\le t\le1)\right)\right.\\
&&\qquad\qquad\qquad\qquad\qquad
\left.-\e\left(\frac1{\gamma h(X_1^\uparrow)}\ind_{\{X_1^\uparrow\ge\eta\}}F(X^\uparrow_t,\,0\le t\le1)\right)\right|\\
&&+\e\left(\frac1{\p({\cal C}_n^{(n)})V^{(n)}(S_n^{(n)\uparrow})}
\ind_{\{S^{(n)\uparrow}_n<\eta\}}F(S^{(n)\uparrow}_{[nt]},\,0\le t\le1)\right)\\
&&\qquad\qquad\qquad\qquad\qquad
+\e\left(\frac1{\gamma h(X_1^\uparrow)}\ind_{\{X_1^\uparrow<\eta\}}F(X^\uparrow_t,\,0\le t\le1)\right)\,.
\end{eqnarray*}
Since $F$ is bounded by a constant, say $B$ and
\begin{equation}\label{alphabeta}
\e\left(\frac1{\p({\cal C}_n^{(n)})V^{(n)}(S_n^{(n)\uparrow})}\right)=1\quad\mbox{and}\quad
\e\left(\frac1{\gamma h(X_1^\uparrow)}\right)=\beta/\gamma\,,\end{equation}
it follows from H\"older's inequality that
the two last terms of the right hand side of the above
inequality are bounded above respectively by $B\p(S^{(n)\uparrow}_n<\eta)$ and $B\p(X_1^\uparrow<\eta)\beta/\gamma$. From
the assumption of convergence and the fact that $\p(X_1^\uparrow>0)=1$, for every $\varepsilon>0$, there exist $n_0$
and  $\eta>0$ such that for all $n\ge n_0$,  $B\p(S^{(n)\uparrow}_n<\eta)<\varepsilon$ and $B\p(X_1^\uparrow<\eta)\beta/\gamma<\varepsilon$.
Finally,  note that from the hypothesis of
convergence and Lemma \ref{harmfunct},  we easily derive that for all $\eta>0$,
\begin{eqnarray*}
&&\e\left(\frac1{\p({\cal C}_n^{(n)})V^{(n)}(S_n^{(n)\uparrow})}
\ind_{\{S^{(n)\uparrow}_n\ge\eta\}}F(S^{(n)\uparrow}_{[nt]},\,0\le t\le1)\right)\\
&&\qquad\qquad\qquad\qquad{\longrightarrow}\;\;
\e\left(\frac1{\gamma h(X_1^\uparrow)}\ind_{\{X_1^\uparrow\ge\eta\}}F(X^\uparrow_t,\,0\le t\le1)\right)\,,
\qquad\mbox{as $n\rightarrow\infty$.}
\end{eqnarray*}
Then we have proved that
\begin{eqnarray*}
&&\e\left(\frac1{\p({\cal C}_n^{(n)})V^{(n)}(S_n^{(n)\uparrow})}F(S^{(n)\uparrow}_{[nt]},\,0\le t\le1)\right)\\
&&\qquad\qquad\qquad\qquad\qquad\qquad{\longrightarrow}\;\;
\e\left(\frac1{\gamma h(X_1^\uparrow)}F(X^\uparrow_t,\,0\le t\le1)\right),
\qquad\mbox{as $n\rightarrow\infty$.}
\end{eqnarray*}
Taking $F\equiv1$ in this relation and comparing with (\ref{alphabeta}), we obtain $\beta=\gamma$, which proves the result.
\end{proof}

\end{document}